\documentclass[12pt]{article}
\usepackage{amssymb,amsthm,amsmath,graphics,psfrag,geometry,longtable}

\newcommand{\Ord}{\mathrm{Ord}}
\newcommand{\SOrd}{\mathrm{SOrd}}
\newcommand{\G}{\mathrm{G}}
\renewcommand{\min}{\mathrm{min}}
\renewcommand{\max}{\mathrm{max}}
\newcommand{\lub}{\mathrm{lub}}
\newcommand{\Pow}{\mathrm{Pow}}
\newcommand{\dom}{\mathrm{dom}}
\newcommand{\ran}{\mathrm{ran}}
\newcommand{\fun}{\mathrm{fun}}
\newcommand{\unique}{\mathrm{unique}}
\newcommand{\chain}{\mathrm{chain}}
\newcommand{\cut}{\mathrm{cut}}
\newcommand{\pred}{\mathrm{pred}}
\newcommand{\ext}{\mathrm{ext}}
\newcommand{\fund}{\mathrm{fund}}
\renewcommand{\top}{\mathrm{top}}
\renewcommand{\bot}{\mathrm{bot}}
\newcommand{\unitop}{\mathrm{unitop}}
\newcommand{\unibot}{\mathrm{unibotsuc}}
\newcommand{\set}{\mathrm{set}}
\newcommand{\Set}{\mathrm{Set}}
\newcommand{\mor}{\mathrm{mor}}
\newcommand{\V}{\mathrm{V}}
\newcommand{\comp}{\mathrm{comp}}
\newcommand{\recmin}{\mathrm{recmin}}
\newcommand{\Gcl}{\mathrm{Gcl}}
\newcommand{\I}{\mathrm{I}}
\newcommand{\N}{\mathrm{N}}
\newcommand{\id}{\mathrm{id}}
\newcommand{\SO}{\mathrm{SO}}
\newcommand{\ZFC}{\mathrm{ZFC}}

\theoremstyle{plain}
\newtheorem{thm}{Theorem}
\newtheorem{lemma}{Lemma}

\theoremstyle{definition}
\newtheorem{rem}{Remark}
\newtheorem{df}{Definition}

\begin{document} 
\author{Peter Koepke, Martin Koerwien}
\title{The Theory of Sets of Ordinals}
\maketitle
\abstract{
We propose a natural theory SO axiomatizing the class of sets 
of ordinals in a model of ZFC set theory. Both theories 
possess equal logical strength. Constructibility theory in SO 
corresponds to a natural recursion theory on ordinals.}

\section{Introduction.}

\footnotetext{2000 \emph{Mathematics Subject Classification;} 03E45}

Cantorian set theory and its axiomatizations describe a 
universe of \emph{hierarchical} sets. According to Cantor's dictum
\begin{quote}
Unter einer ``Menge'' verstehen wir jede Zusammenfassung $M$ von
bestimmten wohlunterschiedenen Objekten $m$ unsrer 
Anschauung oder
unseres Denkens (welche die ``Elemente'' von $M$ genannt werden) zu
einem Ganzen.
\footnote{By a ``set'' we understand any collection into a whole
$M$ of definite and separate objects $m$ of our intuition or our
  thought. These objects are called the ``elements'' of $M$.}
\cite{cantor}
\end{quote}
a set can have (other) sets as its elements and thus one is led to 
the consideration of sets of sets, sets of sets of sets, and so 
on. Such hierarchical sets allow the formalization of the 
fundamental notions of set theory and mathematics:
Kuratowski \cite{kuratowski} defines the ordered 
pair $(x,y)$ as $\{\{x\},\{x,y\}\}$;
and von Neumann \cite{neumann} builds up the ordinal numbers as $0=\emptyset$, 
$1=\{\emptyset\}$, $2=\{0,1\} = \{\emptyset,\{\emptyset\}\}$, 
$3=\{0,1,2\}=\{\emptyset,\{\emptyset\},\{\emptyset,\{\emptyset\}\}\}$, etc.

Whereas the Cantorian notion of set allows to formalize all of 
mathematics in the small language $\{\in\}$, there are some
drawbacks. For example, in axiomatic set theory a consequence of the 
hierarchical notion of set is the familiar but very involved 
recursive definition of the forcing relation for atomic formulae:
$p\Vdash^* \dot{x}=\dot{y}$ (see e.g. 
\cite[Chapter VII,{\S} 3]{kunen}).

It is well-known that a model of Zermelo-Fraenkel set 
theory with the axiom of choice is determined by its sets of 
ordinals \cite[Theorem 13.28]{jech2002}. Also, most constructions in 
set theory can be reduced to constructions of ``flat'' 
sets of previously existing objects. This motivates the present 
article:

In Chapter 2, we define and study a natural 
theory of sets of ordinals 
(SO) which is as strong as the system ZFC and which can 
serve as a foundation of mathematics in a way similar to ZFC. 
The theory SO is two-sorted: {\it ordinals} are taken as given 
atomic objects, avoiding von Neumann's hierarchical ordinals. 
The second sort corresponds to {\it sets of ordinals}. The 
fundamental notion of {\it pairing} is present in the form of G{\"o}del's 
ordinal pairing function 
corresponding to the canonical well-ordering of
$\Ord\times\Ord$ (see \cite[Section 3]{jech2002}).
In Chapter 3 we give definitions for an SO-model within a ZFC-model and
for a ZFC-model within an SO-model. These operations are inverse to each other
and show that ZFC and SO possess the same axiomatic strength.

It is interesting to transfer parts of
standard axiomatic set theory to SO.
In Chapter 4, we carry
out {\it constructibility theory} within SO using
a specific kind of recursion theory on ordinals ($*$-{\it recursion}) which might be
of independent interest.

This article employes a range of canonical coding techniques. In
the interest of space and time we concentrate our exposition upon central
ideas and problems and leave out a great number of technical
details.
The results of this paper were obtained as part of the second 
author's masters thesis \cite{koerwien}, 
supervised by the first author. The 
preparation of the article was financially supported by the 
Mathematical Institute of the University of Bonn.

\section{The theory SO}

Let $L_\SO$ be the language 
$$L_{\SO}:=\{\Ord,\SOrd,<,=,\in,\G\}$$
where $\Ord$ and $\SOrd$ are 
unary predicates, $<$, $=$ and $\in$ are binary
predicates and $\G$ is a two-place function. 
To simplify notation, we use lower case 
greek letters
to range over elements of Ord and lower case roman 
letters to range over elements
of SOrd, so, e.g., 
$\forall\alpha\phi$ stands for $\forall\alpha(\Ord(\alpha)\to\phi)$. Let 
$\alpha\leq\beta$ abbreviate the expression $\alpha<\beta\lor\alpha=\beta$ and $\exists^{=1}$ postulate
the existence of a unique object. For a formula $\phi$, the notation
$\phi(X_1,\dots,X_n)$ means that the set of free variables of $\phi$ is a subset of
$\{X_1,\dots,X_n\}$. SO is the theory axiomatized by the following set of axioms:
\begin{quote}
        \begin{itemize}
        \item [(SOR)] Axiom of sorts\\
        $\forall X,Y((\Ord(X)\leftrightarrow\lnot \SOrd(X))\land\\
	(X<Y\to \Ord(X)\land \Ord(Y))\land\\
        (X\in Y\to\Ord(X)\land \SOrd(Y))\land \Ord(\G(X,Y)))$
        \item [(WO)] Well-ordering axiom\\
        $\forall\alpha,\beta,\gamma(\lnot\alpha<\alpha\land(\alpha<\beta\land\beta<\gamma\to
        \alpha<\gamma)\land\\
	(\alpha<\beta\lor\alpha=\beta\lor\beta<\alpha))\land\\
        \forall a(\exists\alpha(\alpha\in a)\to\exists\alpha(\alpha\in a\land
        \forall\beta(\beta<\alpha\to\lnot\beta\in a)))$
        \item [(INF)] Axiom of infinity (existence of a limit ordinal)\\
        $\exists\alpha(\exists\beta(\beta<\alpha)\land\forall\beta
        (\beta<\alpha\to\exists\gamma(\beta<\gamma\land\gamma<\alpha)))$
        \item [(EXT)] Axiom of extensionality\\
        $\forall a,b(\forall\alpha(\alpha\in a\leftrightarrow\alpha\in b)\to a=b)$
        \item [(INI)] Initial segment axiom\\
        $\forall\alpha\exists a\forall\beta(\beta<\alpha\leftrightarrow\beta\in a)$
        \item [(BOU)] Boundedness axiom\\
        $\forall a\exists\alpha\forall\beta(\beta\in a\to\beta<\alpha)$
        \item [(GPF)] Pairing axiom (G{\"o}del Pairing Function)\\
        $\forall\alpha,\beta,\gamma(\G(\beta,\gamma)\leq\alpha\leftrightarrow$
        $\forall\delta,\epsilon((\delta,\epsilon)<^*(\beta,\gamma)$
        $\to\\ \G(\delta,\epsilon)<\alpha))$\\[2mm]
        Here $(\alpha,\beta)<^*(\gamma,\delta)$ stands for\\
        $\exists\eta,\theta(\eta=\max(\alpha,\beta)\land\theta=
	\max(\gamma,\delta)\land(\eta<\theta\lor\\
        (\eta=\theta\land \alpha<\gamma)
	\lor(\eta=\theta\land \alpha=\gamma\land \beta<\delta)))$,\\
        where $\gamma=\max(\alpha,\beta)$ abbreviates 
        $(\alpha>\beta\land\gamma=\alpha)\lor(\alpha\leq\beta\land\gamma=\beta)$
        \item [(SUR)] $\G$ is onto\\
        $\forall\alpha\exists\beta,\gamma(\alpha=\G(\beta,\gamma))$
        \item [(SEP)] Axiom schema of separation: For all $L_{SO}$-formulae 
        $\phi(\alpha,P_1,\dots,P_n)$ postulate:\\
        $\forall P_1,\dots,P_n\forall a\exists b\forall\alpha(\alpha\in b\leftrightarrow\alpha
        \in a\land\phi(\alpha,P_1,\dots,P_n))$
        \item [(REP)] Axiom schema of replacement: For all $L_{SO}$-formulae
        $\phi(\alpha,\beta,P_1,\dots,P_n)$ postulate:\\
        $\forall P_1,\dots,P_n(\forall\xi,\zeta_1,\zeta_2(\phi(\xi,\zeta_1,P_1,\dots,P_n)\land\\
        \phi(\xi,\zeta_2,P_1,\dots,P_n)\to\zeta_1=\zeta_2)\to\\
	\forall a\exists b\forall\zeta
        (\zeta\in b\leftrightarrow\exists\xi\in a\ \phi(\xi,\zeta,P_1,\dots,P_n)))$
        \item [(POW)] Power set axiom\\
        $\forall a\exists b(\forall z(\exists\alpha(\alpha\in z)\land
        \forall\alpha(\alpha\in z\to\alpha\in a)\to\\
        \exists^{=1}\xi\forall\beta(\beta\in z\leftrightarrow \G(\beta,\xi)\in b)))$\\
        \end{itemize}
\end{quote}
  It is obvious that the structure composed of the ordinals and sets of ordinals
  in ZFC basically satisfies SO (for technical details see proposition
  \ref{prop:SOinZFC}).
  Note that the power set axiom
  of SO postulates the existence of \emph{well-ordered} power sets and thus
  also captures in a certain way the axiom of choice.

We list some 
observations and conventions. Assume SO for the rest
 of this chapter. \emph{$\alpha$
 is an ordinal} and \emph{$a$ is a set} will mean that $\Ord(\alpha)$
 and $\SOrd(a)$ respectively.
We will make use of the class term notation $A=\{x|\phi(x)\}$ familiar
from standard set theory to denote classes of ordinals and sets.
If $A=\{x|\phi(x)\}$ is a non-empty \emph{class of ordinals}, i.e. $\forall x\in
 A(\Ord(x))$, let $\min(A)$ denote the minimal element of $A$. The
 existence of such an element follows from the
 axioms (INI), (SEP) and (WO).
(BOU) ensures the existence of an upper bound for each set $a$, the least of which will
 be noted $\lub(a)$. By (INI) the classes
 $\iota_\alpha:=\{\beta|\beta<\alpha\}$ are sets. 
Using (SEP) and (INI), one sees that the union and 
intersection of two sets are again sets. Finite sets are denoted by
 $\{\alpha_0,\alpha_1,\dots,\alpha_{n-1}\}$. Their existence is
 implied by (INI) and (SEP).
We write $\Pow(b,a)$ for $b$ being a set satisfying
(POW) for $a$. $\xi_b$ will then be the function which to each nonempty subset 
$z$ of $a$ assigns the unique ordinal number $\xi_b(z)$ such that
$\forall\beta(\beta\in z\leftrightarrow \G(\beta,\xi_b(z))\in b)$. $\omega$ denotes
the least element of the class of
limit numbers which by (INF) is not empty. Finally let 
$0:=\min(\{\alpha|\Ord(\alpha)\})$,
$1:=\lub(\{0\})$, etc.

The inverse functions $\G_1$, $\G_2$ of $\G$ are defined via the properties 
$\alpha=\G_1(\beta):=\exists\gamma(\beta=\G(\alpha,\gamma))$ resp.
$\alpha=\G_2(\beta):=\exists\gamma(\beta=\G(\gamma,\alpha))$.
The axioms (GPF) and (SUR) imply the well-known properties of the G{\"o}del pairing function and
its projections, such as bijectivity and monotonicity properties.
To simplify notation, let $(\alpha,\beta):=\G(\alpha,\beta)$. Every set can be regarded as a set of pairs
$a=\{(\alpha,\beta)|(\alpha,\beta)\in a\}$ or more general as a set of $n$-tuples.
In this way $n$-ary relations and functions on ordinals can be encoded as sets.

\begin{df}
  Let $X$, $Y$, $f$, $g$ be sets or classes.
        \begin{eqnarray*}
        \Ord & := & \{\alpha|\Ord(\alpha)\}\\
        \SOrd & := & \{\alpha|\SOrd(\alpha)\}\\
        S & := & \Ord\cup\SOrd\\
        \emptyset & := & \iota_0\\
        \dom(X) & := & \{\alpha|\exists\beta((\alpha,\beta)\in X)\}\\
        \ran(X) & := & \{\beta|\exists\alpha((\alpha,\beta)\in X)\}\\
        \widehat{X} & := & \dom(X)\cup\ran(X)\\
        \fun(X) & := & \forall\alpha,\beta_1,\beta_2((\alpha,\beta_1)\in X\land(\alpha,\beta_2)\in X
        \to\beta_1=\beta_2)\\
        f:X\to Y & := & \fun(f)\land\dom(f)=X\land\ran(f)\subset Y\\
        \alpha=f(\beta) & := & (\alpha,\beta)\in f\\
        g\circ f & := & \{(\alpha,\beta)|\exists\gamma(\gamma=f(\alpha)\land\beta=g(\gamma))\}\\
        \alpha~Y~\beta & := & (\alpha,\beta)\in Y\\
        X\times Y & := & \{\gamma|\G_1(\gamma)\in X\land\G_2(\gamma)\in Y\}\\
        X\restriction Y & := & \{(\alpha,\beta)\in X|\alpha\in Y\}\\
        f''X & := & \{\beta|\exists\alpha\in X\ (\alpha,\beta)\in f\}
        \end{eqnarray*}
\end{df}

\begin{thm}[Transfinite induction]
        Let $\phi(\alpha,X_1,\dots,X_n)$ be an $L_{SO}$-formula. Then for all $X_1,\dots,X_n$,
        $$\forall\alpha(\forall\beta<\alpha(\phi(\beta,X_1,\dots,X_n))\to\phi(\alpha,X_1,\dots,X_n))$$
        implies
        $$\forall\alpha\phi(\alpha,X_1,\dots,X_n)$$
\end{thm}
\begin{proof}
        Otherwise, by (WO), there would be a minimal counterexample $\alpha$ contradicting the assumption.
\end{proof}

\begin{thm}[Transfinite recursion]
        Let $R:\Ord\times\SOrd\to\Ord$ be a function defined by some
        formula $\phi(\alpha,f,\beta,X_1,\dots,X_n)$. 
        Then there exists a unique function $F:\Ord\to\Ord$ defined
        by a formula $\psi(\alpha,\beta,X_1,\dots,X_n)$ such that
        \begin{equation}
          \label{eq1}
          \forall\alpha(F(\alpha)=R(\alpha,F\restriction\iota_\alpha))
        \end{equation}
\end{thm}
\begin{proof}
        This is proved similar to the recursion theorem in ZF:
        We define the notion of \emph{approximation functions} which are
        set-functions defined on proper initial segments of Ord, satisfying 
        (\ref{eq1}) on their domain. Then we obtain $F$ as the union of all of these
        approximation functions.
\end{proof}

\begin{rem}
As in ZF this result can be generalized from the relation $<$ to arbitrary set-like well-founded relations.
\end{rem}
To give an example how to work inside SO and what kind of problems can arise, we define the structure
of real numbers with addition and multiplication. This procedure indicates the potential of SO to serve as a foundational theory of mathematics, similar to ZFC.

Using the recursion theorem, we define addition and multiplication on ordinals.
This provides us with
the structure $(\mathbb{N},+,\cdot,0,1)$,
$\mathbb{N}:=\iota_\omega$,
which satisfies the axioms of second order Peano Arithmetic.

The standard construction of the rational numbers by equivalence
classes of tuples runs into trouble because these equivalence classes are
sets and cannot be assembled together to be
the \emph{set} of rationals. One solves
this problem by representing the equivalence classes by their minimal elements.
Then addition and multiplication can be defined on these representatives in the obvious
way.

To define the real numbers as Dedekind cuts of the rational numbers, we take,
by the power set axiom 
(POW) a set $y$ such that $\Pow(y,\mathbb{Q})$.
Via the function $\xi_y$ we can assign to each (non-empty) subset of $\mathbb{Q}$ 
a unique ordinal. Let $\mathbb{R}$
be the set of all such ordinals whose corresponding subset
of $\mathbb{Q}$ is the left half of a Dedekind cut. Then addition and multiplication
on $\mathbb{R}$ can be defined as in the usual theory of Dedekind
cuts.

Standard structures and constructions
such as topological spaces,
Cartesian products, quotient spaces are available in SO.
The formation of sets of sets can usually be avoided by representing 
equivalence classes by minimal representatives. 
Some constructions, however,
are no longer canonical due to the non-uniqueness of 
power sets in SO.

\section{The bi-interpretability of SO and ZFC}

We introduce a syntactical notion of {\it inner model} for arbitrary 
first order
languages and of interpretations of formulae in those inner models. That notion
is contained as a special case in the definition of 
{\it interpretability} as introduced
in \cite{hodges}. 

\begin{df}
  \label{Def:InnerModels}
        Let $L_1$ and $L_2$ be first order languages and $T_1$ an
        $L_1$-theory. 
        Let $I$ be the (index-)set consisting of the non-logical symbols of
        $L_2$ (including the identity relation) together with another symbol $u$. A collection
        $\mathcal{A}=(\phi_i)_{i\in I}$ of $L_1$-formulae 
        is called a \emph{$T_1$-definable $L_2$-structure} if
        \begin{itemize}
        \item[(i)] $\phi_u$ has exactly one free variable $x$ and
          $T_1\vdash\exists x\ \phi_u(x)$. We write $x\in U$ instead of $\phi_u(x)$.
        \item[(ii)] For all relation symbols $r\in I$ 
          the free variables of $\phi_r$ are exactly $v_1,\dots,v_n$ where $n$ is
          the arity of $r$.
        \item[(iii)] For all function symbols $f\in I$ the free 
          variables of $\phi_f$ are exactly $v_1,\dots,v_{n+1}$ where $n$ is
          the arity of $r$. Moreover, 
          $T_1\vdash\forall v_1,\dots,v_{n+2}\in
          U(\phi_f(v_1,\dots,v_n,v_{n+1})\land
          \phi_f(v_1,\dots,v_n,v_{n+2})\to\phi_=(v_{n+1},v_{n+2}))$
          and\\
          $T_1\vdash\forall v_1,\dots,v_n\in
          U\exists v_{n+1}\in U\ \phi_f(v_1,\dots,v_n,v_{n+1})$.
        \item[(iv)] For all constant symbols $c\in I$, $\phi_c$ has exactly
          one free variable $x$ and
          $T_1\vdash\exists x\in U(\phi_c(x))\land\forall x,y\in U(\phi_c(x)\land \phi_c(y)\to\phi_=(x,y))$.
        \item[(v)] $T_1$ proves that $\phi_=$ defines a
          \emph{congruence relation}
          for $L_2$, i.e. it has the properties of
          an equivalence relation and respects all functions and relations defined by the formulas of $\mathcal{A}$.
        \end{itemize}
\end{df}

\begin{df}
        Let $L_1$ and $L_2$ be first order languages, 
$T_1$ an $L_1$-theory and $\mathcal{A}$ a 
        $T_1$-definable $L_2$-structure.
        Then for an $L_2$-formula $\psi$ the 
        \emph{relativization  of $\psi$ to $\mathcal{A}$}
        is an $L_1$-formula $\psi^\mathcal{A}$ defined 
        by recursion on the structure of $\psi$:
        \begin{itemize}
          \item[(i)] If $\psi\equiv(x=y)$, where $x$ and $y$ are
          variables, then $\psi^\mathcal{A}:=\phi_=(x,y)$.
          \item[(ii)] If $x$ is a variable, $c$ is a constant symbol and $\psi\equiv(x=c)$ then 
            $\psi^\mathcal{A}:=\phi_c(x)$.
          \item[(iii)] If $x$ is a variable, $f$ is a function symbol, 
            $t_1,\dots,t_n$ are $L_2$-terms
            and $\psi\equiv(x=f(t_1,\dots,t_n))$ 
then
          $\psi^\mathcal{A}:=\\
	\exists x_1,\dots,x_n\in U((x_1=t_1)^\mathcal{A}
      \land\dots\land(x_n=t_n)^\mathcal{A}\land\phi_f(x_1,\dots,x_n,x))$.
          \item[(iv)] 
If $r$ is a relation symbol (including the identity) then \\
          $r(t_1,\dots,t_n)^\mathcal{A}:=
            \exists x_1,\dots,x_n\in U((x_1=t_1)^\mathcal{A}\land\\
	\dots\land(x_n=t_n)^\mathcal{A}\land\phi_r(x_1,\dots,x_n))$
          \item[(v)] $(\lnot\psi)^\mathcal{A}:=\lnot\psi^\mathcal{A}$, 
            $(\psi_1\lor\psi_2)^\mathcal{A}:=\psi_1^\mathcal{A}\lor\psi_2^\mathcal{A}$ and
            $(\exists x\psi)^\mathcal{A}:=\\
	\exists x\in U(\psi^\mathcal{A})$.
        \end{itemize}
        If $\Phi$ is a set of $L_2$-formulae we define $\Phi^\mathcal{A}:=\{\phi^\mathcal{A}|\phi\in\Phi\}$.
\end{df}

\begin{df}[Interpretability]

        Let $L_1$ and $L_2$ be first order languages, $T_1$ an $L_1$-theory and $T_2$ an $L_2$-theory. Then $T_2$ is 
        \emph{interpretable in $T_1$} (or \emph{$T_1$ interprets $T_2$}) iff 
        there is a $T_1$-definable $L_2$-structure $\mathcal{A}$
        such that $T_1\vdash T_2^\mathcal{A}$.

        $T_1$ and $T_2$ are \emph{bi-interpretable} iff $T_1$ interprets $T_2$ and $T_2$ interprets $T_1$.
\end{df}

\begin{rem}
        If $T_2$ is interpretable in $T_1$ and $T_1$ is consistent
        then $T_2$ is consistent.
\end{rem}


\begin{thm}
  \label{prop:SOinZFC}
        ZFC interprets SO.
\end{thm}
\begin{proof}
        The SO-ordinals will be interpreted by the ordinals in ZFC. To distinguish the set of ordinals 
        $\{\alpha|\alpha<\beta\}$ from the ordinal $\beta$, we interpret the SO-sets of ordinals by the class
        $\SOrd:=\{x\cup\{\Omega\}|x\subset\Ord\}$, i.e., we ``mark''
        the sets of ordinals by a fixed set $\Omega$ 
which is not an ordinal, e.g., $\Omega:=\{\{\emptyset\}\}$.
        
The relations and functions of $L_{SO}$ can
        be defined on $${\mathrm{S}}(V):=\Ord\cup\SOrd$$ in the obvious way.
        Clearly the theory SO is designed to describes the properties of 
        ordinals and sets of ordinals in
a ZFC-model,
        so the validity of the axioms is immediately verified. Note
        that the proof of $(POW)^{{\mathrm{S}}(V)}$ requires the axiom of choice since we
        obtain a power set according to SO 
from a well-ordering of the
        corresponding ZFC-power set.
\end{proof}

We claim that ZFC and SO are bi-interpretable. 
So we have 
to define a model $V(S)$ of ZFC in a given SO-model $S$.
First we motivate our construction.

Given a set $a$ in a ZFC-universe, the structure $A:=(\mathrm{TC}(\{a\}),\in\restriction A\times A)$
of its transitive closure determines uniquely this set.
This structure has some obvious properties: it is well-founded, extensional,
has a unique minimal element (the empty set) and a unique top
element $a$ such that for all other elements $b\in\mathrm{TC}(\{a\})$ there
exists a descending $\in$-chain $c=(c_0,\dots,c_n)$ from $a$ to $b$ 
such that $b=c_n\in c_{n-1}\in\dots\in c_0=a$.

From now on we will work in SO.
As was remarked above, SO-sets can be regarded as sets of pairs, i.e., as binary relations.
The class of all binary relations 
satisfying the properties of the last paragraph 
will be the universe
of our model $\V(S)$. 
We shall define appropriate identity and element relations on $V(S)$.

\begin{df}
	$ $
        \begin{longtable}{lll}
        $\pred_a(\alpha)$ & $:=$ & $\{\beta|\beta~a~\alpha\}$\\
        $\ext(a)$ & $:=$ & $\forall\alpha,\beta\in\widehat{a}(\pred_a(\alpha)=\pred_a(\beta)\to\alpha=\beta)$\\
        $\fund(a)$ & $:=$ & $\forall b(b\subset\widehat{a}\land 
        b\neq\emptyset\to\exists\beta\in b\forall\alpha\in b(\lnot\alpha~a~\beta))$\\
        $\unique(c,a)$ & $:=$ & $\forall\alpha,\beta,\gamma\in 
        c((\alpha~a~\gamma\land \beta~a~\gamma)\lor$\\
	& & $\qquad(\gamma~a~\alpha\land \gamma~a~\beta)\to\alpha=\beta)$\\
        $\chain(c,a,\alpha,\beta)$ & $:=$ & $\unique(c,a)\land\alpha\in c
        \land\beta\in c\land$\\
	& & $\qquad\forall\gamma\in c(\lnot\gamma=\alpha\leftrightarrow\exists\delta
	\in c(\gamma~a~\delta))\land$\\
        & & $\qquad\forall\gamma\in c(\lnot\gamma=\beta\leftrightarrow\exists\delta
	\in c(\delta~a~\gamma))$\\ 
        $\alpha=\top(a)$ & $:=$ & $\forall\beta\in\widehat{a}\exists c(\chain(c,a,\alpha,\beta))$\\
        $\unitop(a)$ & $:=$ & $\exists\alpha(\alpha=\top(a))$\\
        $\alpha=\bot(a)$ & $:=$ & $\pred_a(\alpha)=\emptyset$\\
        $\unibot(a)$ & $:=$ & $\forall\alpha(\alpha=\bot(a)\to\exists^{=1}\beta(\alpha~a~\beta))$\\
        \pagebreak $\set(a)$ & $:=$ & $a\neq\emptyset\land\fund(a)\land\ext(a)\land$\\
	& & $\qquad\unitop(a)\land\unibot(a)$\\
        $\Set$ & $:=$ & $\{a|\set(a)\}$
        \end{longtable}
\end{df}

        If $a\neq\emptyset$, $\fund(a)$ implies the existence of an 
        $\alpha$ such that $\alpha=\bot(a)$. If we have $\ext(a)$, $\bot(a)$ is uniquely
        defined. Also if $\top(a)$ exists, $\fund(a)$ implies that it
        must be unique.

        Many elements of $\Set$ correspond to the same transitive closure of
        a set. We have to define an appropriate equivalence relation $\approx$ on $\Set$.

\begin{df}
        \begin{eqnarray*}
        \mor(f,a,b) & := & f:\widehat{a}\to\widehat{b}\land\forall\alpha,\beta\in\widehat{a}(\alpha~a~\beta\to 
        f(\alpha)~b~f(\beta))\\
        id_a & := & \{(\alpha,\alpha)|\alpha\in a\}\\
        a\approx_{f,g}b & := & \mor(f,a,b)\land\mor(g,b,a)\land g\circ 
        f=id_{\widehat{a}}\land f\circ g=id_{\widehat{b}}\\
        a\approx b & := & \exists f,g(a\approx_{f,g}b)
        \end{eqnarray*}
\end{df}

\begin{df}
        \begin{eqnarray*}
        \cut(a,\alpha) & := & \{\beta\in a|\exists
        c_1,c_2(\chain(c_1,a,\alpha,\G_1(\beta))
        \land\\
	& & \qquad\chain(c_2,a,\alpha,\G_2(\beta)))\}\\
        a~\widetilde{\in}~b & := &
        \set(b)\land\exists\beta~b~\top(b)(a\approx\cut(b,\beta))\\
        \V(S) & := & (\Set,\approx,\widetilde{\in})
        \end{eqnarray*}
\end{df}

\begin{thm}
  \label{SOintZFC}
  SO interprets ZFC.
\end{thm}
To prove the Theorem, it is enough to show $\SO\vdash \ZFC^{\V(S)}$.
As a first observation, $\approx$ is a congruence relation
for $\widetilde{\in}$. An easy induction shows that
the validity of $L_{\in}$-formulas interpreted in $\V(S)$ does not depend on the choice of the
representatives of the parameters.

A code for the empty set exists, namely $\widetilde{\emptyset}:=\{(0,1)\}\in\Set$: 
If $a~\widetilde{\in}~\widetilde{\emptyset}$ then there must be 
an $\alpha~\widetilde{\emptyset}~\top(\widetilde{\emptyset})$ such
that $a\approx\cut(\widetilde{\emptyset},\alpha)$. 
Then obviously $\alpha=0$ and $a\approx\emptyset$ which implies $a=\emptyset$ and
$\lnot\set(a)$ by definition.

\begin{rem}
  We have decided to represent sets by codes of relations that do not mention
  explicitely their carriers. Since this does not allow to distinguish
  between the empty carrier and one-element carriers, we represent the
  empty set by the two-element carrier relation
  $\widetilde{\emptyset}=\{(0,1)\}$.
  Since every non-empty transitive set has the empty set as its
  $\in$-minimal element, we find a copy of $\widetilde{\emptyset}$ at
  the bottom of every representative. This is expressed by the unibotsuc-condition.
\end{rem}

Next we prove a lemma by which we can define elements of Set with
prescribed $\widetilde{\in}$-predecessors.

\begin{lemma}
        \label{lemma:ausw}
        Let $a$, $d$ be sets such that $\fund(a)$, $\ext(a)$ 
        and $\unibot(a)$ and let $\emptyset\neq d\subset\widehat{a}\setminus\{\bot(a)\}$.
        Then for arbitrary $\alpha\in\Ord\setminus\widehat{a}$ 
        $$\set(a,d,\alpha):=\{(\delta,\alpha)|\delta\in d\}\cup\bigcup\limits_{\delta\in d}\cut(a,\delta)$$
        is an element of $\Set$ and for all $b\in\Set$ we have
        $$b~\widetilde{\in}~\set(a,d,\alpha)\Leftrightarrow\exists\delta\in d(b\approx\cut(a,\delta))$$
\end{lemma}
\begin{proof}
        We have $\widehat{\set(a,d,\alpha)}=\{\alpha\}\cup\bigcup\limits_{\delta\in d}\widehat{\cut(a,\delta)}\subset
        \{\alpha\}\cup\widehat{a}$ and for 
        $\beta,\gamma\in\bigcup\limits_{\delta\in d}\widehat{\cut(a,\delta)}$ obviously 
        \begin{equation}
          \label{eq2}
          \beta~a~\gamma\text{ iff }\beta~\set(a,d,\alpha)~\gamma\text{.}
        \end{equation}
        The property $\fund(\set(a,d,\alpha))$ is clear
        because if $b\cap\widehat{a}\neq\emptyset$ for 
        $\emptyset\neq b\subset\widehat{\set(a,d,\alpha)}$ then an $a$-minimal
        element of $b\cap\widehat{a}$ is $\set(a,d,\alpha)$-minimal in $b$.

        To prove $\ext(\set(a,d,\alpha))$, first observe that for 
        $\beta\in\widehat{\set(a,d,\alpha)}\setminus\{\alpha\}$
        $\pred_a(\beta)=\pred_{\set(a,d,\alpha)}(\beta)$ since there
        is a $\delta\in d$ such that $\beta\in\widehat{\cut(a,\delta)}$
        and $\pred_a(\beta)\subset\widehat{\cut(a,\delta)}$. So if we have $\beta,\gamma$ with
        $\pred_{\set(a,d,\alpha)}(\beta)=\pred_{\set(a,d,\alpha)}(\gamma)$ 
        then the case where $\beta,\gamma\in\widehat{a}$ is trivial
        because of $\ext(a)$. There remains the case where 
        $\beta=\alpha$, $\gamma\in\widehat{a}$. Then there exists $\delta\in d$ such that
        $\gamma\in\widehat{\cut(a,\delta)}$ and since 
        $\delta\in\pred_{\set(a,d,\alpha)}(\beta)=\pred_{\set(a,d,\alpha)}(\gamma)$ also 
        $\delta~\set(a,d,\alpha)~\gamma$ and so a chain $c$ from 
        $\delta$ down to $\gamma$ would have no $\set(a,d,\alpha)$-minimal
        element, contradicting $\fund(\set(a,d,\alpha))$.

        As for $\unitop(\set(a,d,\alpha))$, first of all we have\\
        $\chain(\{\alpha\},\set(a,d,\alpha),\alpha,\alpha)$.
        If then $\beta\in\widehat{\set(a,d,\alpha)}\setminus\{\alpha\}$, there 
        exists $\delta\in d$ and $c$ such that 
        $\chain(c,a,\delta,\beta)$ and thus $\chain(c,\set(a,d,\alpha),\delta,\beta)$ by (\ref{eq2}). Then
        obviously $\chain(c\cup\{\alpha\},\set(a,d,\alpha),\alpha,\beta)$.

        It is clear that $\bot(\set(a,d,\alpha))=\bot(a)$. As\\
        $\set(a,d,\alpha)\setminus\{(\delta,\alpha)|\delta\in d\}\subset a$
        and $\lnot(\bot(\set(a,d,\alpha))~\set(a,d,\alpha)~\alpha)$ 
        (since $\bot(\set(a,d,\alpha))\notin d$), we have $\unibot(\set(a,d,\alpha))$.

        To show
        $b~\widetilde{\in}~\set(a,d,\alpha)\Leftrightarrow\exists\delta\in d(b\approx\cut(a,\delta))$, 
        it suffices to prove that for $\delta\in d$, $\cut(a,\delta)=\cut(\set(a,d,\alpha),\delta)$.
        But this is clear by (\ref{eq2}).
\end{proof}

Instead of 
formal proofs of the ZFC axioms relativised to $\V(S)$, 
we just indicate the main ideas; many details are
routine and trivial.

To prove the 
\emph{scheme of separation} let
$\phi(b,x_1,\dots,x_n)$ 
be an $L_\in$-formula and $a,x_1,\dots,x_n\in\Set$.
We put
$$d:=\{\gamma\in\pred_a(\top(a))|\phi^{\V(S)}(\cut(a,\gamma),x_1,\dots,x_n)\}$$ 
which is a set by (SEP). If $d=\emptyset$ or $a=\widetilde{\emptyset}$ 
then clearly $\widetilde{\emptyset}$ is the set we are looking for. 
Otherwise $d$ satisfies the conditions of the preceding lemma
and $\set(a,d,\alpha)$ for some $\alpha\notin\widehat{a}$ has the desired properties.

The proofs of the \emph{axioms of choice} and \emph{union} are similar. If in
$\V(S)$, $a$ is a set of non-empty pairwise disjoint sets, we obtain a
choice-set by applying Lemma \ref{lemma:ausw} to 
$d:=\{\gamma|\exists\beta~a~\top(a)(\gamma=\min(\pred_a(\gamma)))\}$. 
As for the union of a set $a$, apply Lemma \ref{lemma:ausw} with
$d:=\bigcup\limits_{\gamma~a~\top(a)}\pred_a(\gamma)$.

For the proof
of the \emph{axiom of extensionality}, consider $a,b\in\Set$ with equal sets of
$\widetilde{\in}$-predecessors. We therefore have already unique 
isomorphisms between the $\cut(a,\alpha)$ and
corresponding $\cut(b,\beta)$ parts of $a$ and $b$ for $\alpha~a~\top(a)$ and $\beta~b~\top(b)$. Taking the union of
all these isomorphisms (noting that they are compatible)
and mapping $\top(a)$ to $\top(b)$ gives the desired isomorphism of $a$ and $b$.

For the \emph{axiom of foundation},
suppose that there were an infinite decending $\widetilde{\in}$-chain
beginning with
$a_0\in\Set$. Then all elements of this chain are represented 
in $a_0$ as $\cut(a_0,\alpha_i)$ for certain $\alpha_i\in\widehat{a_0}$ which
results in an infinite descending $a_0$-chain contradicting
$\fund(a_0)$.
 
The \emph{axiom of infinity} can be proved by explicitly
constructing a code of the (ZFC-)ordinal number $\omega$ as 
$$\widetilde{\omega}=
\{(0,1)\}\cup\{(\alpha,\beta)|0<\alpha<\beta\leq\omega\}.$$
This set exists by (INF), (INI) and (SEP) and clearly
$\widetilde{\omega}\in\Set$.

As \emph{pairing} follows from replacement and infinity, only replacement 
and the power set axiom remain to be shown.

For the construction of the \emph{power set} $p$ of a set $a$, we first take
a power set $b$ of $\pred_a(\top(a))$ by (POW) and obtain
a numbering of the subsets of $\pred_a(\top(a))$ by the $\xi_b$
function. We avoid the possible complication that these numbers could be elements of $\widehat{a}$
by replacing them by their images under the bijection 
$\alpha\mapsto(\alpha,\zeta)$ for some fixed
$\zeta\geq\lub(\widehat{a})$. Fix a new top element $\gamma$ (an
arbitrary number not colliding with any number that appears in our
construction). We would like to take the union
of all $\cut(a,\alpha)$ for $\alpha~a~\top(a)$, of all 
$(\alpha,(\xi_b(z),\zeta))$ for $\emptyset\neq z\subset\pred_a(\top(a))$
and $\alpha\in z$ and of all $((\xi_b(z),\zeta),\gamma)$.
But this union possibly does not satisfy extensionality because there could be  $\alpha\in\widehat{a}$ 
such that $z:=\pred_a(\alpha)\subset\pred_a(\top(a))$. 
Thus, in these cases we have to replace $(\xi(z),\zeta)$ by these $\alpha$. 
Finally we have to take care of the fact that the empty set belongs 
to the power set by adding $(\beta,\gamma)$ to our relation, where $\beta$ is the unique
successor of $\bot(a)$. In this way we obtain a set which satisfies the defining property
of the power set of $a$.

\emph{Replacement} is the most involved schema to prove. 
Given a formula $\phi(b_1,b_2,X_1,\dots,X_n)$ such that $\phi^{\V(S)}$
is functional and a set $a$,
we have to ``unify'' all $b_2$ such that 
$\phi^{\V(S)}(b_1,b_2,X_1,\dots,X_n)$ for
$b_1~\widetilde{\in}~a$. But these $b_2$ are only determined
up to isomorphism. So we have to find uniform representatives 
for these sets. Using (ERS) we can show that there exists an
$\alpha$ such that all relevant $b_2$ are represented by sets 
$b_2^\prime$ such that $\widehat{b_2^\prime}\subset\iota_\alpha$ and thus all
these relations are subsets of $\iota_\alpha\times\iota_\alpha$. 
Then we take a power set $y$ for $\iota_\alpha$, thus enumerating
all $b_2^\prime$. With respect to this order, we cobble the 
$b_2^\prime$ (which we make disjoint by the method described in the preceding paragraph)
together to form one relation, i.e., at each step we add the part not 
yet represented and the links to what is already constructed.
Then we obtain a set $b$ which satisfies the requirements of the lemma. We put
$$d:=\{\alpha\in\widehat{b}|\exists
x~\widetilde{\in}~a(\phi^{\V(S)}(x,\cut(b,\alpha),X_1,\dots,X_n))\}$$
and apply Lemma \ref{lemma:ausw} to find the desired set. This
completes the proof of Theorem \ref{SOintZFC}.

\begin{thm}
        Assuming ZFC, there exist an isomorphism
        $F_{VS}:\V({\mathrm{S}}(V))\simeq V$. Assuming SO, there exists an
        isomorphism $F_{SV}:{\mathrm{S}}(\V(S))\simeq S$.
\end{thm}
\begin{rem}
  The notion of ``isomorphism'' in this 
theorem has to be
  understood in the following way:

  The three parts ``being a function'', ``being one-one'' and ``being
  onto'' must be formulated in the appropriate language. The
  statement of the theorem is that these formulae are consequences of the
  corresponding theory. For example, ``$F_{VS}$ is a function''
  translates to ZFC $\vdash\forall x,y\in\V({\mathrm{S}}(V))(x\approx
  y\to F_{VS}(x)=F_{VS}(y))$, where the function symbol
  $F_{VS}$ must be replaced by
a $\in$-formula that definies this function.
\end{rem}
\begin{proof}
        Working in a ZFC-model $V$ we remark that the elements of 
        $\Set$ constructed in ${\mathrm{S}}(V)$ are, seen as ZFC relations, extensional
        and well-founded. So they can be collapsed uniquely to 
        transitive sets (before collapsing we remove the bot-element). The top-element
        of this transitive set will be defined to be the image of 
        an application $F_{VS}$ and $F_{VS}$ is easily seen to be an isomorphism.

        Now starting from an SO-model $S$ we can, as in the proof of 
        the ZFC infinity axiom, define canonical representatives for
        the ``ordinal numbers'' in $\V(S)$ and thus also canonical 
        representatives for the sets of ordinals. In that way we can obtain
        ${\mathrm{S}}(\V(S))$ as an SO-class with SO-definable 
        relations and functions. Then we define $F_{SV}$ by assigning to a code
        for an ordinal the rank (after removing the bot-element) 
        of its top-element and to a code for a set of ordinals the set 
        of images of its $\widetilde{\in}$-elements. Again the 
        proof that this defines an isomorphism is straightforward.
\end{proof}

\section{$*$-recursion and the constructible\\ model $S^*$}

In this chapter, we sketch how to carry out {\it constructibility theory}
in the framework of SO. We present a notion of $*$-{\it recursiveness}
which generalizes the ordinary recursive functions from $\omega$ to
Ord. We shall see that the $*$-recursive sets of ordinals are exactly
the constructible sets of ordinals.
 
\begin{df}
A function from a cartesian product of Ord into Ord is 
\emph{$*$-recursive} if it is generated by the following schema
        \begin{itemize}
        \item[(i)] For all $m\leq n<\omega$ the following functions are $*$-recursive:\\[2mm]
        \begin{tabular}{ll}
        $\id:\Ord\to\Ord$, & $\beta\mapsto\beta$,\\
        $\pi_m^n:\Ord^n\to\Ord$, & $(\beta_0,\dots,\beta_{n-1})\mapsto\beta_m$,\\
        $f_\lor:\Ord^2\to\Ord$, & $(\beta,\gamma)\mapsto\left\{\begin{array}{cl}
        1,\text{ if} & \beta>0\lor\gamma>0\\0,\text{ if} & \beta=0\land\gamma=0\end{array}\right.$,\\
        $f_<:\Ord^2\to\Ord$, & $(\beta,\gamma)\mapsto\left\{\begin{array}{cl}
        1,\text{ if} & \beta<\gamma\\0,\text{ if} & \beta\geq\gamma\end{array}\right.$,\\
        $f_=:\Ord^2\to\Ord$, & $(\beta,\gamma)\mapsto\left\{\begin{array}{cl}
        1,\text{ if} & \beta=\gamma\\0,\text{ if} & \beta\neq\gamma\end{array}\right.$,\\
        $f_\lnot:\Ord\to\Ord$, & $\beta\mapsto\left\{\begin{array}{cl}
        1,\text{ if} & \beta=0\\0,\text{ if} & \beta>0\end{array}\right.$,\\
        $G_1:\Ord\to\Ord$, & $\beta\mapsto G_1(\beta)$\\
        $G_2:\Ord\to\Ord$, & $\beta\mapsto G_2(\beta)$\\
        $G:\Ord^3\to\Ord$, & $(\beta,\gamma,\delta)\mapsto\left\{\begin{array}{cl} 
        \G(\gamma,\delta),\text{ if} & \G(\gamma,\delta)<\beta\\ \beta,\text{ if} & \G(\gamma,\delta)\geq\beta\end{array}\right.$
        \end{tabular} 
        \item[(ii)]
        Let $g:\Ord^n\to\Ord$ and $h_i:\Ord^m\to\Ord$ ($i\in\{1,\dots,n\}$) be $*$-recursive.
        Then the \emph{composition}\\
        $\comp(g,h_1,\dots,h_n):\Ord^m\to\Ord$,\\
        $(\beta_1,\dots,\beta_m)\mapsto g(h_1(\beta_1,\dots,\beta_m),\dots,h_n(\beta_1,\dots,\beta_m))$\\
        is $*$-recursive.
        \item[(iii)]
        Let $g:\Ord^{n+m}\to\Ord$ and $h_i^j:\Ord^n\to\Ord$ ($i\in\{1,\dots,n-1\}$, $j\in\{1,\dots m\}$) be
        $*$-recursive. Then the \emph{recursive minimization} of $g$ (w.r.t. the $h_i^j$) is $*$-recursive:\\
        $f:=\recmin(g,(h_1^1,\dots,h_{n-1}^1),\dots,(h_1^m,\dots,h_{n-1}^m)):\Ord^n\to\Ord$,\\
        $(\beta_0,\dots,\beta_{n-1})\mapsto\left\{\begin{array}{cl}
        \min(\{\delta<\beta_0|g(\delta,\vec{\beta},f^{(\beta_0,\dots,\beta_{n-1})}(\beta_0,\vec{\gamma}_1),
        \dots,\\ \quad
        f^{(\beta_0,\dots,\beta_{n-1})}(\beta_0,\vec{\gamma}_m))>0\})\text{,
        if defined}&\\ \\ \beta_0
        \text{, otherwise}\end{array}\right.$\\
        where 
        $\vec\beta:=\beta_1,\dots,\beta_{n-1}$, 
        $\vec\gamma_j:=h_1^j(\delta,\vec\beta),\dots,h_{n-1}^j(\delta,\vec\beta)$
        and\\
        $f^{(\beta_0,\dots,\beta_{n-1})}(\alpha_0,\dots,\alpha_{n-1})$ is defined as\\
        $~~~~~~~~~~~~~~\left\{\begin{array}{cl}
        f(\alpha_0,\dots,\alpha_{n-1})\text{, if
        }\alpha_0\leq\beta_0,\dots,\alpha_{n-1}\leq\beta_{n-1}\\ \quad
        \text{ and
        }(\alpha_0,\dots,\alpha_{n-1})\neq(\beta_0,\dots,\beta_{n-1})&\\ \\
        0\text{, otherwise}\end{array}\right.$
        \end{itemize}
\end{df}
        \begin{rem}
          \label{RemDefInSO}
        As one may expect, $*$-recursion can be formally defined in
        the theory SO. First we define a \emph{reasonable} numbering
        of the functions using the G{\"o}del pairing function,
        that is, one that allows to recover the inductive definition of a
        $*$-recursive function by the projections $\G_1$ and $\G_2$. 
        Then we can define recursively a function $ari$ which yields
        the arity of the function coded by an ordinal number, and 
        a function $\mathit{FUN}(\alpha,\beta)$
        which assigns to a code of a function $\alpha$ and an ordinal 
        number $\beta$ the value of the coded function at the argument
        $\beta$ (regarded as a tuple 
        $\beta=(\beta_1,(\beta_2,(\dots,(\beta_{n-1},\beta_n)\dots)))$ for $n$-ary functions).

        Here some technical difficulties arise, as the arity of the 
        functions can increase during the recursive computation
        if the function is defined by recursive minimization or composition. By the properties of the
        G{\"o}del pairing function, this can cause an increase of the
        argument $\beta$. 

        This problem can be solved either by
        defining $\mathit{FUN}$ by recursion on the well-founded 
        relation $<^*$ defined by $\gamma<^*\delta:=
        (\G_1(\gamma)=\G_1(\delta)\land G_2(\gamma)<\G_2(\delta))
        \lor(\G_1(\gamma)<\G_1(\delta)\land\G_2(\gamma)<\Gcl(\G_2(\delta)))$,
        where $\eta:=\Gcl(\G_2(\delta))$ is the minimal ordinal number
        greater than or equal to $G_2(\delta)$ which is closed under $\G$,
        i.e. $\forall\eta_1,\eta_2<\eta(\G(\eta_1,\eta_2)<\eta)$. 

        Another solution of the problem is to restrict $\mathit{FUN}$ to
        arguments $\beta<\beta_0$ where $\beta_0$ is greater than all 
        arguments needed in the present context (which is sufficient
        for our purposes). Then the numbering of the $*$-recursive 
        functions can be defined such that the codes always dominate the 
        arguments when using composition or recursive minimization 
        (modify the codes by something like $\text{new code}=
        \G(\beta_0,\text{old code})$). Since the codes become
        smaller during the recursive computation and dominate 
        the arguments, the decrease of the arguments of 
        $\mathit{FUN}$ is guaranteed by the properties of the 
        G{\"o}del pairing function and the computation works.
        \end{rem}

        \begin{df}
        Let $\mathit{rec}$ be the class of ordinal codes for $*$-recursive 
        functions (in the sense of the preceding remark). Let
        $\alpha,\beta,\gamma,\delta\in\Ord$ and $a\in\SOrd$. Then define

        \begin{eqnarray*}
                \I(\alpha,\beta,\gamma) & := &
                \{\eta<\alpha|\beta\in\mathit{rec}
                \land\mathit{FUN}(\beta,(\eta,\gamma))>0\}\\
                \SOrd^*(a) & := & \exists\delta(a=\I(\delta))\\
                \SOrd^* & := & \{a|\SOrd^*(a)\}\\
                \N(a) & := & \min({\delta|a=\I(\delta)})\\
                \N & := & \{\delta|\exists a\in\SOrd^*(\delta=\N(a))\}
        \end{eqnarray*}
        \end{df}

        $\SOrd^*$ is the class of $*$-recursively definable 
        (in short $*$-definable) sets. $\N$ is the class of (minimal) names for
        $*$-definable sets.
        In the above definition, $\gamma$ plays the role of a parameter
        (or a tuple of parameters, using the G{\"o}del pairing function).
        \begin{df}
          We say that a class $\SOrd^\prime\subset\SOrd$ defines 
          an \emph{inner model} of SO if $S^\prime:=\Ord\cup\SOrd^\prime$
          satisfies SO under the obvious interpretation 
          (here we use the symbol $S^\prime$ to denote the $L_{SO}$-substructure
          with domain $S^\prime$).
        \end{df}

        \begin{thm}
          \label{SO*InnerModel}
          $\SOrd^*$ defines an inner model which we denote by $S^*$.
        \end{thm}
        We sketch roughly the main arguments for the proof of Theorem \ref{SO*InnerModel}.
        First of all, one observes that many of the axioms of SO only 
        concern ordinal numbers and thus are absolute for all
        inner models. Also the proof of (INI) is trivial. As (SEP)
        follows easily from (REP),
        the only axioms that need proof are (POW) and (REP).

        The following fact is crucial for the proofs of (POW) and
        (REP):\\[2mm]{\bf Fact: }
        {\it
          The notion of $*$-recursion can be defined
          $*$-recursively, i.e., the functions $\mathit{ari}$ and $\mathit{FUN}$
          are definable as $*$-recursive functions. More precisely,
          there is a \emph{universal} $*$-recursive function
          $\mathit{FUN}$ such that for any $*$-recursive $f$ there is
          an $\alpha$ such that
          $$f(x_0,\dots,x_{n-1})=\mathit{FUN}(\alpha,(x_0,\dots,x_{n-1})).$$
        }
        
\noindent In fact, the schema of recursive minimization is built exactly in a way to make this
        possible (separated schemas of minimization and recursion as
        in ordinary recursion theory seem
        not to be sufficient).

        We very briefly note some techniques used
        for the $*$-recursive definition of $\mathit{FUN}$. First, by the
        projections $\pi_m^n$ and composition, the arity and the order of arguments
        of every $*$-recursive function can be modified arbitrarily. 
        By the functions $f_\lor$, $f_\lnot$, $f_<$, $f_=$ and composition, ``conditions''
        can be formulated $*$-recursively. Then clearly definitions 
        by cases are possible using recursive minimization like
        ``the minimal number $\delta$ such that
        ($\mathit{condition}_1>0$ 
        and $\delta=\mathit{value_1}$) or 
        ($\mathit{condition}_2>0$ and $\delta=\mathit{value_2}$)
        etc.''. Now we can define $\mathit{FUN}$ as a recursive minimization
        of a function which distinguishes the different cases 
        (atomic functions, composition, recursive minimization). An important
        point is that $*$-recursion can deal uniformly with tuples of 
        arbitrary length (e.g. argument tuples), 
        treating them as sequences by recursively defined projection 
        functions which yield the n-th component of a tuple
        (the important difference to the $\pi_m^n$ functions is, 
        that $m$ and $n$ become arguments of the function).

        We return to the proof of Theorem \ref{SO*InnerModel}.
        In order to prove (POW), take a $*$-recursively definable set
        $a$. We have to find a $*$-recursive function which defines
        a \emph{power set} for $a$. For that, we shall be able 
        to test $*$-recursively if an ordinal defines a subset of $a$. 
        Using the $*$-recursive version of $\mathit{FUN}$, we can define a
        function $e(\alpha,\beta)$ which returns $0$ or $1$ depending
        on whether $\alpha$ is an element of the $*$-recursive set
        defined by $\beta$ (regarded as a triple) or not. Then define
        $$s(\alpha,\beta,\eta):=(\mathit{min}~\epsilon<\eta(e(\epsilon,\alpha)\land 
        f_\lnot(e\epsilon,\beta)))=\eta$$
        using recursive minimization ($\eta$ shall be an arbitrary sufficiently large number and the symbol
        $\land$ should be replaced by applications of $f_\lor$ and $f_\lnot$).
        This function tests if $I(\alpha)\subset I(\beta)$.
        Now we are able to define $*$-recursively a predicate that
        expresses that two numbers define the same $*$-recursive
        set which allows us to express that an 
        ordinal number is a minimal name for a $*$-definable set.
        Finally we can define the desired power set by a function 
        $g(\alpha,\beta)$ which is defined to return $1$ if $\beta$ is
        a minimal name of a subset of $a$ and $\alpha$ is an element 
        of this subset, and which otherwise returns $0$ 
        (of course, formally the set $a$ has to be expressed by a
        name, which becomes a parameter in the definition).
        The function $g$ can be defined in SO using the second
        approach descibed in Remark \ref{RemDefInSO} since the class
        of minimal names for $a$ can be bounded using (POW) and (REP).

        For the proof of (REP) let $\phi(\alpha,\beta,x_1,\dots,x_n)$ 
        be an $L_{SO}$-formula with parameters in $S^*$ such
        that $\phi^{S^*}$ is functional. Let $a\in\SOrd^*$. The set 
        $\{\beta|\exists\alpha\in
        a(\phi^{S^*}(\alpha,\beta,x_1,\dots,x_n)\}$ 
        must be shown to be $*$-definable. Since $\phi$ is an
        arbitrary formula, we have to find a way to $*$-recursively 
        calculate the truth of formulas (that are relativized to
        $S^*$). Before we continue the proof of (REP),
        we state a theorem that corresponds to the well-known
        reflection principle in ZFC. Its proof is similiar to the
        proof of \cite[Theorem 12.14]{jech2002}
        \begin{thm}
          \label{Reflection}
                Let $\psi(\alpha_1,\dots,\alpha_n,b_1,\dots,b_m)$ 
                be an $L_{SO}$-formula and let $\phi:=\psi^{S^*}$.
                Then there exists an ordinal number $\alpha$ such that
                $$\forall\alpha_1,\dots,\alpha_n\in\iota_\alpha\forall
                b_1,\dots,b_m\in\{d\in\SOrd^*|d\subset\iota_\alpha\}
                (\phi\leftrightarrow\phi^\alpha)$$
                where $\phi^\alpha$ is recursively defined as 
                the formula $\phi$ with all quantifiers restricted to $\iota_\alpha$
                resp. $\{d\in\SOrd^*|d\subset\iota_\alpha\}$.
        \end{thm}
        We choose a reasonable numbering of all
        formulas including constants for ordinal numbers and elements
        of $\SOrd^*$ (represented by their names). To distinguish 
        ordinal numbers from sets, we fix a maximal height $\alpha$
        for ordinals we want to deal with and code sets by 
        $(\text{name},\alpha)$. Then we define $*$-recursively 
        a function $\mathit{subst}(\beta,n,\gamma)$ (realized 
        as a recursive minimization of a definition by cases) that substitutes 
        the variable $v_n$ in the formula $\beta$ (i.e. $\beta$ is a number of a formula)
        by the constant $\alpha$. Now a function that calculates the
        truth of formulas can be defined (again a recursive minimization of a definition by cases)
        which in the quantifier-case substitutes the quantified 
        variable by a constant that makes the formula true if
        this is possible. The bound $\alpha$ can be found by Theorem
        \ref{Reflection} applied to $\phi^{S^*}$. This concludes the proof of
        Theorem \ref{SO*InnerModel}.
        
        $S^*$ is not only an inner model of SO, it is the smallest inner model,
        i.e. for all inner models $S^\prime\subset S$
        we have $S^*\subset S^\prime$. This can be seen quite easily by
        the absoluteness of the definition of $\mathit{FUN}$ which implies
        the absoluteness of all classes $\I(\delta)$, hence of $\SOrd^*$
        which therefore must be included in all inner models.

        \begin{rem}
                In the following, inner models of ZFC are always understood to be definable, transitive
                and to contain all ordinal numbers.
        \end{rem}

        We conclude this paper by showing that $S^*$ corresponds to the constructible universe $L$ of
        ZFC.

        \begin{thm}
                \label{thm:inn}
                Let $V$ be a model of ZFC and $S$ be a model of SO.
                \begin{itemize}
                \item[(i)] If $M\subset V$ is an inner model then 
                  ${\mathrm{S}}(M)\subset{\mathrm{S}}(V)$ is an inner model.
                \item[(ii)] If $S^\prime\subset S$ is an inner model 
                  then $\V(S^\prime)\subset\V(S)$ is an inner model.
                \end{itemize}
        \end{thm}
        \begin{proof}
                The inclusion ${\mathrm{S}}(M)\subset{\mathrm{S}}(V)$ is evident.
                $\V(S^\prime)\subset\V(S)$ follows from the absoluteness
                of $\Set$, $\approx$ and $\widetilde{\in}$.
                
                We just prove (i), the proof of (ii) being quite similar.
                By assumption we have $\mathrm{ZFC}^M$
                and by inner interpretability $(\mathrm{SO}^{{\mathrm{S}}(V)})^M$,
                i.e., $\mathrm{SO}^{{\mathrm{S}}(M)}$. Also ${\mathrm{S}}(M)$ 
                an \emph{in ${\mathrm{S}}(V)$} definable class, since it can be
                shown that ${\mathrm{S}}(M)=F_{SV}^{{\mathrm{S}}(V)}[{\mathrm{S}}(M^{\V({\mathrm{S}}(V))})]$,
                where $F_{SV}^{{\mathrm{S}}(V)}$ denotes the isomorphism $F_{SV}$ defined in ${\mathrm{S}}(V)$ and 
                $M^{\V({\mathrm{S}}(V))}:=\\
	\{x\in\V({\mathrm{S}}(V))|\phi^{\V({\mathrm{S}}(V))}(x)\}$
                if $M=\{x|\phi(x)\}$.
        \end{proof}
        
        With our methods we can prove a version of \cite[Lemma 13.28]{jech2002}.

        \begin{thm}
                \label{thm:incl}
                Let $V$ be a model of ZFC and $M,N\subset V$ inner models.
                If $\{x\in M|x\subset\Ord\}\subset\{x\in N|x\subset\Ord\}$,
                then $M\subset N$. As a corollary, if $M$ and $N$
                have the same sets of ordinals, they are identical.
        \end{thm}
        \begin{proof}
                Obviously, ${\mathrm{S}}(M)\subset{\mathrm{S}}(N)$ and by the absoluteness of
                $\Set$ also $\V({\mathrm{S}}(M))\subset\V({\mathrm{S}}(N))$. Since the isomorphism
                $F_{VS}:\V({\mathrm{S}}(V))\to V$ can be shown to be absolute, the valid assertion 
                $(\forall x\exists y(x=F_{VS}(y)))^M$ implies
                $\forall x\in M\exists y\in M(x=F_{VS}(y))$. Now let
                $x\in M$ and $y$ such that $x=F_{VS}(y)$. This implies
                $y\in\V(S(M))$ and consequently $y\in\V(S(N))$.
                The valid assertion $(\forall y\in\V({\mathrm{S}}(V))\exists\bar{x}(\bar{x}=F_{VS}(y)))^N$
                now implies
                $\forall y\in\V({\mathrm{S}}(N))\exists\bar{x}\in N(\bar{x}=F_{VS}(y))$
                and finally $x=F_{VS}(y)=\bar{x}\in N$.
        \end{proof}

        By the last two results above, $*$-recursive sets are exactly the constructible sets of ordinals:

        \begin{thm}
                \label{thm:AAA}
                Let $S$ be a model of SO and let $L^{\V(S)}$ denote
                the constructible inner model of $\V(S)$. 
                Then $S^*=F_{SV}[{\mathrm{S}}(L^{\V(S)})]$.
        \end{thm}
        \begin{proof}
                $S^*\subset F_{SV}[{\mathrm{S}}(L^{\V(S)})]$ is clear, as
                $F_{SV}[{\mathrm{S}}(L^{\V(S)})]$ is easily seen to be an inner model of $S$.
                By minimality of the constructible universe and theorem
                \ref{thm:inn} we have $L^{\V(S)}\subset\V(S^*)$.
                Then ${\mathrm{S}}(L^{\V(S)})\subset{\mathrm{S}}(\V(S^*))$ and
                $F_{SV}[{\mathrm{S}}(L^{\V(S)})]\subset F_{SV}[{\mathrm{S}}(\V(S^*))]=S^*$.
        \end{proof}

        \begin{thm}
                Let $S$ be a model of SO. Then $\V(S^*)=L^{\V(S)}$.
        \end{thm}
        \begin{proof}
                Theorem \ref{thm:AAA} yields
                ${\mathrm{S}}(\V(S^*))={F_{SV}}^{-1}[S^*]=\\
	{F_{SV}}^{-1}[F_{SV}[{\mathrm{S}}(L^{\V(S)})]]={\mathrm{S}}(L^{\V(S)})$.
                This implies that $\V(S^*)$ and $L^{\V(S)}$ have
                the same sets of ordinals. Now apply theorem
                \ref{thm:incl}.
        \end{proof}

\bibliographystyle{alpha}
\bibliography{Article1.bib}

\begin{thebibliography}{Kun80}

\bibitem[Can95]{cantor}
Georg Cantor.
\newblock Beitr\"age zur Begr\"undung der transfiniten Mengenlehre.
\newblock {\em Mathematische Annalen}, 46:481--512, 1895.

\bibitem[Hod93]{hodges}
Wilfrid Hodges.
\newblock {\em Model Theory}.
\newblock Cambridge University Press, 1993.

\bibitem[Jec02]{jech2002}
Thomas Jech.
\newblock {\em Set Theory, Third Edition}.
\newblock Springer, 2002.

\bibitem[Koe01]{koerwien}
Martin Koerwien.
\newblock {Die Theorie der Ordinalzahlmengen und ihre Beziehung zur
  G{\"o}delschen Konstruktibilit{\"a}tstheorie}.
\newblock Diplom thesis, Universit{\"a}t Bonn, 2001.

\bibitem[Kun80]{kunen}
K.~Kunen.
\newblock {\em Set Theory. An Introduction to Independence Proofs}.
\newblock North Holland, 1980.

\bibitem[Kur21]{kuratowski}
Casimir Kuratowski.
\newblock Sur la notion de l'ordre dans la théorie des ensembles.
\newblock {\em Fundamenta Mathematicae}, 2, 1921.

\bibitem[vN61]{neumann}
John von Neumann.
\newblock {\em Collected Works, Vol. I}.
\newblock Pergamon Press, 1961.

\end{thebibliography}

\bigskip

\begin{tabular}{lr}
Peter Koepke & Martin Koerwien\\
Mathematisches Institut & Equipe de Logique Math\'ematique\\
Universit\"at Bonn     & UFR de Math\'ematiques (case 7012)\\
Beringstra{\ss}e 1  & Universit\'e Denis Diderot Paris 7\\
D-53115 Bonn & 2 place Jussieu\\
email: koepke@math.uni-bonn.de & F-75251 Paris Cedex 05\\
 & email: koerwien@logique.jussieu.fr
\end{tabular}

\end{document}